\documentclass[12pt]{amsart}
\usepackage{amsmath,amsfonts,amsthm,amsopn,cite,mathrsfs}
\usepackage{epsfig,verbatim}
\usepackage{subfigure}

\newcommand{\stft}{short-time Fourier transform}

\newcommand{\bdl}{band-limited}

\newcommand{\rep}{representation}

\newtheorem{tm}{Theorem}%[section]
\newtheorem{lemma}[tm]{Lemma}
\newtheorem{prop}[tm]{Proposition}

\newcommand{\rem}{\noindent\textsl{REMARK:}}

 \theoremstyle{definition}
 \newtheorem{definition}{Definition}

\newcommand{\beqa}{\begin{eqnarray*}}
\newcommand{\eeqa}{\end{eqnarray*}}

\DeclareMathOperator*{\supp}{supp}

\newcommand{\field}[1]{\mathbb{#1}}
\newcommand{\bR}{\field{R}}        %  real numbers
 \def\cH{\mathcal{H}}
 
 \def\cE{\mathcal{E}}

 \def\cX{\mathcal{X}}

\def\rd{\bR^d}

\def\lrd{L^2(\rd)}

\def\<{\left<}
\def\>{\right>}

\def\inv{^{-1}}

\def\mv1{M_v^1}

\newcommand{\fsh}{f^\sharp}
\newcommand{\card}{\mathrm{card}\,}
\newcommand{\vf}{\varphi}
\newcommand{\cY}{\mathcal{Y}}

\begin{document}
\begin{abstract}
We show that the homogeneous approximation property and the comparison
theorem hold for arbitrary coherent
frames. This observation answers some questions about the density of
frames that are not covered by the theory of Balan, Casazza, Heil, and
Landau. The proofs are a variation of the method developed by
Ramanathan and Steger. 
\end{abstract}

\title[Homogeneous Approximation Property and Coherent Frames]{The Homogeneous Approximation Property and the Comparison
  Theorem  for Coherent Frames}
\author{Karlheinz Gr\"ochenig}
\address{Faculty of Mathematics \\
University of Vienna \\
Nordbergstrasse 15 \\
A-1090 Vienna, Austria}
\email{karlheinz.groechenig@univie.ac.at}
\date{}
\thanks{K.~G.~was supported by the Marie-Curie Excellence Grant MEXT-CT 2004-517154}
\maketitle

%\section{Introduction}

Frames provide redundant non-orthogonal expansions in Hilbert space,
intuitively they should therefore be ``denser'' than orthonormal
bases. The first fundamental result is  Landau's necessary condition
for sets of sampling of \bdl\ 
functions~\cite{landau67}. A second   fundamental idea is  the
homogeneous 
approximation property and the comparison theorem of Ramanathan and
Steger~\cite{RS95} in the context of Gabor frames. Since then  many
contributions have investigated the density of frames and
varied and applied  the method of Ramanathan and
Steger. See~\cite{CDH99,FNP07,GR96,sch98} for a sample of papers. The
approach in~\cite{RS95}  
culminates in the deep investigation of Balan, Casazza, Heil, and
Landau~\cite{BCHL06,BCHL06a} who have found a density theory for a general class of
frames that are labeled by a discrete Abelian group. 

However, the BCHL theory is not universally applicable, for instance,
the density of wavelet frames does not fall under this general theory,
and a separate investigation is necessary, as was shown in a sequence
of papers  by Heil and
Kutyniok~\cite{HK03,HK07} and by Sun~\cite{Sun07,SunZ04,Sun07b}. Their work uses quite
explicit and long  calculations in the group of affine
transformations.

In  this note we  show that the  approach of  Ramanathan and Steger
yields the homogeneous approximation 
property  and a comparison theorem for the class of coherent
frames. These frames  arise as
subsets of the orbit of a square-integrable group representation.
 In particular,  wavelet frames possess this structure. The advantage
 of the abstract point of view is the conciseness and simplicity of
 the proofs and the much more general range of validity. The 
 additional insight is that the homogeneous approximation
 property and the comparison theorem are a consequence of structure
 alone,
 we do not need to assume additional localization properties as in
 \cite{BCHL06a}.

We will work in the context of a locally compact group $G$  with Haar
measure $dx$.  We 
write $|U| = \int _G \chi _U(x) \, dx $ for the measure of a set
$U\subseteq G$ and $K^c$ for the complement of $K$ in $G$. 

\emph{Amalgam spaces.} 
Let $U=U\inv \subseteq G $ be a symmetric compact neighborhood of the
identity element $e\in G $ and let 
$$
f^\sharp (x) = \sup _{y\in xU} |f(y)|
$$
be the local maximum function of $f$. 
Then  the left amalgam
space  $W(L^\infty , L^2)$ is defined by the norm 
\begin{equation}
  \label{eq:1}
\|f \|_{W(L^\infty , L^2)} =   \|f^\sharp \|_2 \, .
\end{equation}

Recall that a set $\cX \subseteq G $  is called relatively (left) separated
if 
\begin{equation}
  \label{eq:3}
  \sup _{x\in G} \mathrm{card}\, (\cX \cap xU) <\infty \, .
\end{equation}
for some, hence all, compact neighborhood(s) $U$ of $e$. 
Equivalently, for any compact set $K\subseteq G$ the sum 
$\sum _{x\in \cX } \chi _{xK} $ is uniformly bounded in $G$. 

\begin{lemma}
  \label{lem:1}
Let $\cX = \{x_j\}$ be a relatively (left) separated set  in $G$ and
$K\subseteq G $ a compact set. Then 
\begin{equation}
  \label{eq:2}
  \sum _{x_j \not \in K} |f(x_j)|^2 \leq C \, \int _{K^c U } |\fsh (x)
  |^2 \, dx \, .
\end{equation}
\end{lemma}

\begin{proof}
  Set $C_0 = \|\sum _{j\in J} \chi _{x_j U } \|_\infty  $.
By \eqref{eq:3} this constant is finite. 
First we note that 
$|f(x_j)| \leq \fsh (x)$ whenever $ x_j \in xU$ or $ x\in x_jU$
(because $U$ is symmetric).
Consequently
$$
|f(x_j)|^2 \leq \frac{1}{|U|}\int _{x_j U} \fsh (x)^2 \, dx \, .
$$
So after summing over $j$, we obtain
\begin{eqnarray*}
  \sum _{x_j \not \in K} |f(x_j)|^2 &\leq &  \frac{1}{|U|} \sum _{x_j
    \not \in K} \int _{x_j U} \fsh (x)^2 \, dx \\
&\leq & \frac{1}{|U|} C_0 \, \int _{K^c U} \fsh (x) ^2 \, dx \, ,
\end{eqnarray*}
because $\sum _{x_j \not \in K } \chi _{x_j U } (x) \leq  C_0 \chi
_{K^cU} (x)$ for all $x\in G$. 
\end{proof}

Next we consider the setup for  coorbit theory and coherent frames.
Let $(\pi , \cH )$ be a  unitary,
\emph{square-integrable} \rep\ of $G$ on a Hilbert space $\cH $. 
For fixed $g\in \cH , g\neq 0$, we set 
$V_gf =  \langle f, \pi (x) g\rangle$; this representation coefficient
is the analog of the 
wavelet transform or \stft\ for the group $G$.

We investigate frames of the form $\{\pi (x_j) g : j\in J\}$,
so-called coherent frames. This means that there exist positive
constants $A,B >0$, the frame bounds, such that
\begin{equation}
  \label{eq:4a}
A\|f\|^2 \leq \sum _{j} |\langle f, \pi (x_j)g\rangle |^2 = \sum _{j}
|V_gf(x_j) |^2 \leq B\|f\|^2   
\end{equation}
for all $f\in \cH $. From \eqref{eq:4a} it is clear that the sampling
of the 
representation coefficients with respect to  $g$ on the set $\cX $ 
must be in $\ell ^2$. This is close, but not equivalent, to saying the
$V_gf \in W(L^\infty , L^2)$. Thus the frame condition motivates the
following definition.

\begin{definition}
  A vector (wavelet) $g\in \cH $ is ``nice'', if the mapping $V_g$
  (the generalized wavelet transform) maps $\cH $ into $W(L^\infty
  ,L^2)$.  
\end{definition}
In most situations the set of nice vectors is dense in $\cH $ and
examples can be easily constructed (see Remark 1). 
In the following we assume  that  $\{\pi (x_j) g : j\in J\}$ is a
frame for $\cH $ and that  $\{ h_j : j\in J\}$ is a dual frame.  We
may  take any dual frame, it need not be the canonical dual. 
We remark that $\{ x_j\}$ must be relatively separated in this case
(otherwise the upper frame bound would not exist). 

The intuition behind coherent frames is that the frame vector $\pi
(x_j)g$ lives near the point $x_j$ in the ``$G$-plane''. Consequently,
a vector of the form $\sum _{x_j \in K} c_j \pi (x_j)g$ should live on
$K$. 

To formalize this intuition, we introduce the following subspaces of
$\cH $. Let  $L\subseteq G$ be a compact subset of $G$  and   $y\in
G$, then $V_{L,y}$ is the  subspace of $\cH $ spanned by the vectors
of the dual frame 
$h_j , x_j \in yL$. Formally, 
\begin{equation}
  \label{eq:4}
  V(yL) = \mathrm{span}\, \{ h_j : x_j \in yL\} \, ,
\end{equation}
and we let $P_{L,y}$ be the orthogonal projection from $\cH $ onto
$V(yL)$. 

The following  homogeneous approximation property is a weak form  of
 localization  of a frame.  The proof for coherent frames follows
 closely the proof for bandlimited functions in~\cite{GR96}.

\begin{prop}[Homogeneous Approximation Property]\label{hap}
Assume  that $g$ is nice and that   $\{\pi (x_j) g : j\in J\}$ is a
frame for $\cH $ and that  $\{ h_j : j\in J\}$ is a dual frame.  Then 
for every $f\in \cH $ and $\epsilon >0$ there exists a compact set
$L\subseteq G$ (depending on $f$ and $\epsilon $) such that 
\begin{equation}
  \label{eq:6}
   \sup _{x\in yK} \| \pi (x)
  f - P_{KL,y} \pi (x) f \| < \epsilon \, 
\end{equation}
holds for  \emph{all} compact sets $K\subseteq G$ and all $y\in G$. 
\end{prop}

\begin{proof}
Since $g$ is nice, the representations coefficient $V_gf $ is in
$W(L^\infty , L^2)$ for all $f\in L^2$. Therefore for any $\delta >0$
there is a compact set $L\subseteq G$, such that on $L^cU$, an  open
neighborhood of the complement of $L$,  the maximal function 
$(V_gf)^\sharp $ is small, precisely, 
\begin{equation}
  \label{eq:42}
  \int _{L^cU} (V_gf)^\sharp (x) \, dx < \delta \, .
\end{equation}

We will use  the frame expansion 

\begin{equation}
  \label{eq:43}
  \pi (x) f = \sum _j \langle \pi (x) f, \pi (x_j ) g \rangle h_j \, .
\end{equation}
Recall the following  property of orthogonal projections: If  $P_{KL,y} h = \sum
_{j: x_j \in yKL} c_j h_j$, then  
\begin{equation}
  \label{eq:7}
  \|h-P_{KL,y}y\| \leq \|h-\sum _{j: x_j \in yKL} d_j h_j \| \, ,
\end{equation}
for any choice of coefficients $d_j$. Choosing  the
special coefficients of the frame expansion~\eqref{eq:43}, we obtain
that 
\begin{eqnarray*}
\|\pi (x) f-P_{KL,y}\pi (x)f\|^2 &\leq & \|\pi (x) f -  \sum _{j: x_j
  \in yKL} \langle \pi (x)f, \pi (x_j )g\rangle    h_j \|^2 \\
&= &  \|\sum _{j: x_j
  \not \in yKL} \langle \pi (x)f, \pi (x_j) g\rangle    h_j \|^2 \\
&\leq & A\inv  \sum _{j: x_j   \not \in yKL} |\langle \pi (x)f, \pi (x_j)
g\rangle |^2 \\
&=& A\inv  \sum _{j: x_j   \not \in yKL} |\langle f, \pi (x \inv x_j) 
g\rangle |^2 \, .
\end{eqnarray*}
 The last inequality follows from the fact that $\{ h_j\} $ is
a frame and obeys the Bessel inequality with constant $A\inv $, where
$A $ is the lower frame bound of $\{\pi (x_j)g\}$ in~\eqref{eq:4a}.

Some gymnastics: Since $x_j    \not \in yKL$, we have $x\inv x_j 
\not \in x\inv  yKL $. Further, $\{x\inv x_j \} $ is relatively
separated, and therefore 
$$
\|\sum _j \chi _{x\inv x_j U} \|_\infty = C_0 < \infty \, ,
$$
where this constant is independent of $x$. 

We apply Lemma~\ref{lem:1} and obtain
\begin{eqnarray*}
  \sum _{j: x\inv x_j   \not \in x\inv yKL} |\langle f, \pi (x \inv x_j)
g\rangle |^2 &=&   \sum _{j: x\inv x_j   \not \in x\inv yKL} |V_gf(x \inv x_j
) |^2 \\
 &\leq & C \int _{(x\inv y KL)^c U} \, |V_gf ^\sharp (x) |^2 \, dx  
\end{eqnarray*}
So for $K\subseteq G $ compact and $x\in yK$ we obtain
\begin{equation}
  \label{eq:cc2}
    \| \pi (x) f - P_{KL,y} \pi (x)
  f\|^2 \leq   A\inv C \sup _{x\inv y =z \in K\inv } \int
  _{(zKL)^c U} V_gf ^\sharp (x)^2 \, dx \, .
\end{equation}
% \begin{eqnarray*}
%    \| \pi (x) f - P_{KL,y} \pi (x)
%   f\|^2 &\leq &  A\inv C \sup _{x\inv y =z \in K\inv } \int
%   _{(zKL)^c U} V_gf ^\sharp (x)^2 \, dx \\
% &\leq &  \int  _{\bigcup _{z\in K\inv } (zKL)^cU}
% V_gf ^\sharp (u)^2 du \, . 
% \end{eqnarray*}
If $K\subseteq G $ is compact and $x\in yK$, then $z=x\inv y \in
K\inv$, therefore  $L\subseteq zKL$ for all $z\in K\inv $ and  $(zKL)^c U
\subseteq L^cU$. Therefore 
$$
\sup _{z\in K\inv} \int _{(zKL)^cU}
V_gf ^\sharp (u)^2 du  \leq \int _{L^cU} (V_g f)^\sharp (u)^2 \, du  \, .
$$
If we now choose $L$ such that $CA\inv \int _{L^cU} (V_g f)^\sharp (u)^2 \,
du < \epsilon $, which is possible by \eqref{eq:42}, then we have
proved that 
$$\sup _{x\in yK} \| \pi (x)
  f - P_{KL,y} \pi (x) f \| < \epsilon \, ,
$$
and the homogeneous approximation property is proved. 
\end{proof}

\rem\ In most situations,  nice vectors $g$ exist in abundance. For
instance,  if $(\pi , \cH)$ is \emph{irreducible} and square-integrable, then
a nice vector  can be constructed as follows: 
Let $g_0 \in \cH $ be admissible, and $k$ be continuous with compact
support.  Then $g=\pi (k)g_0= \int k(x) \pi (x)g_0 \, dx$ has
the property that $V_gf ^\sharp \in W(C, L^2)$ for all $f\in \cH
$. This follows from a convolution equation in~\cite{fg89jfa}. See
also~\cite{FG07}. 

If $\pi $ is the reducible  representation by translations and
dilations on $\lrd $ defined by  $\pi
(x,s)f(t) = |\det e^{-sA/2}| f(e^{-sA} (t-x)$,  where $x,t\in \rd,
s\in \bR ^+$ and $A $ is a real-valued $d\times d$-matrix with all
eigenvalues having positive real part, then every Schwartz function
$g$ such that $\supp \hat{g} \subseteq \{ \omega : 0 < a \leq |\omega
| \leq b\}$ is nice. This follows immediately from the fact that the
wavelet transform $V_gg(x,s) = \langle g, \pi (x,s)g\rangle $ has
compact support in $s$.  See~\cite{GKT92} for more on admissible
vectors in this case. 

\vspace{3 mm}

\textbf{The Comparison Lemma.}
Before we estimate the dimension of the finite-dimensional subspaces
defined in~\eqref{eq:4}, we state a simple estimate for the trace of an
operator. 
\begin{lemma}\label{trace}
  Let $T$ be a positive trace-class operator on a Hilbert space $\cH $
  and  $\{h_k: k\in J\}$ be a  frame with frame bounds
  $A, B>0$. Then 
$$
\frac{1}{B}\, \sum _{k} \langle T h_k , h_k\rangle \leq  \mathrm{tr}\, T
\leq \frac{1}{A} \sum _{k} \langle T h_k , h_k\rangle \, .
$$
\end{lemma}

\begin{proof}
  By the spectral theorem there exists an orthonormal set $\vf _j $
  and positive  eigenvalues $\lambda _j$ such that $Tf = \sum _j \lambda _j
\langle f, \vf _j \rangle \vf _j$. Then 
\begin{eqnarray*}
  \sum _{k} \langle T h_k , h_k\rangle &=& \sum _j \sum _k \lambda _j
  |\langle \vf _j , h_k \rangle |^2 \\
&\leq& B \sum _j \lambda _j \|\vf _j \|_2 \\
&=& B \sum _j \lambda _j = B \, \mathrm{tr}\, T \, .
\end{eqnarray*}
Likewise, by using the lower frame bound, we obtain $\sum _{k} \langle
T h_k , h_k\rangle \geq A \, \mathrm{tr}\, T$. 
\end{proof}

Let $\cE _g= \{\pi (x_j ) g : x_j \in \cX \}$ and $\cE _h= \{\pi (y_k
) h : y_k \in \cY \}$ be two coherent frames. As above,  denote the dual frame of
$\cE_g$ by $\{ h_j \}$ and its frame bounds by $A,B$.

For compact subsets $K,L \subseteq G $ and $y\in G $, define the
subspaces 
\begin{eqnarray}
V_g(yKL) &=& \mathrm{span}\, \{ h_j: x_j \in yKL\} \\  
W_h(yK) &=& \mathrm{span}\, \{ \pi (y_j) h : y_k \in yK\} \, ,  
\end{eqnarray}
and let $P=P_g(yKL)$ and $Q=Q_h(yK)$ be the orthogonal projections
onto the subspaces $V_g$ and $W_h$. 
Finally let $T=QPQ: W_h \to W_h $. $T$ is a positive operator of
finite rank. We estimate its trace from above and below. 

To keep track of these subspaces and projections, keep in mind that
$\cE _g$ is a \emph{given frame} whose density we want to understand and that
$\cE _h $ is the \emph{reference frame}. If possible, $\cE _h$ is
chosen to be  a
Riesz basis or an orthonormal basis. This has been  the case in most
applications of the comparison theorem so far.

\begin{tm}
  \label{comp}
Given $\epsilon >0$, there exists  a compact set
$L\subseteq G$, such that for all $y\in G$ and all compact sets $K
\subseteq G$ 
\begin{equation}
  \label{eq:68}
\|h\|^2  B\inv (1-\epsilon)   \card \{k:  y_k \in yK\} \leq \card \{x_j \in yKL\} \, .
\end{equation}
\end{tm}

\begin{proof}
  Since $T$ is a product of orthogonal projections, its eigenvalues
  are between $0$ and $1$ and 
  \begin{equation}
    \label{eq:65}
    \mathrm{tr}\, T \leq \mathrm{rank}\, P \leq \card \{x_j \in yKL\}\, .
  \end{equation}
To obtain a lower bound, we use  Lemma~\ref{trace} and obtain
\begin{eqnarray*}
\mathrm{tr}\, T &\geq & \frac{1}{B} \sum _{k} \langle T \pi
(y_k)h, \pi (y_k)h \rangle  \\
&\geq & \frac{1}{B}  \sum _{k: y_k \in yK} \langle T \pi (y_k)h, \pi
(y_k)h \rangle\\
 &=&   \sum
    _{k: y_k \in yK} \langle QPQ \pi (y_k)h, \pi (y_k)h \rangle \\
&=& \frac{1}{B}  \sum _{k: y_k \in yK} \langle P \pi (y_k)h, \pi (y_k)h \rangle \\
&=& \frac{1}{B}  \sum _{k: y_k \in yK} \Big(\langle  \pi (y_k)h, \pi (y_k)h \rangle +
  \langle (\mathrm{I}-P) \pi (y_k)h, \pi (y_k)h \rangle \\
&=&  \frac{\|h\|^2}{B} \card \{k: y_k \in yK\}  + (*)  \, .
\end{eqnarray*}
To estimate $(*)$, we apply the homogeneous approximation
property (Proposition~\ref{hap}). Choose a compact set $L\subseteq G $ such that 
$$
\|\pi (x)h - P_{KL,y}\pi (x)h\| < \epsilon \|h\|
$$
for all compact sets $K\subseteq G$ and  $x\in yK$. Then 
\begin{eqnarray*}
(*) &\leq & \frac{1}{B}\sum _{k: y_k \in yK} \|\pi (y_k)h - P (\pi
(y_k)h) \| \|\pi (y_k)h\|  \\
&\leq & \frac{1}{B}\sum _{k: y_k \in yK} \epsilon \|h\| \, \|\pi (y_k)h \|
=  \frac{1}{B}\card \{k: y_k \in yK\}   \epsilon \|h\|^2 \, .
\end{eqnarray*}
By combining these estimates we obtain a lower bound
\begin{equation}
  \label{eq:67}
  \sum _{k: y_k \in yK} \langle T \pi (y_k)h, \pi (y_k)h \rangle  \geq  \card \{k:
  y_k \in yK\} \|h\|^2  B\inv (1-\epsilon) \, .
\end{equation}
Combining \eqref{eq:65} and \eqref{eq:67}, the claim is proved. 
\end{proof}

% \begin{cor}
%   A coherent frame cannot have arbitrarily large gaps. Specifically,
%   if $\{ \pi (x) g: x\in \cX\}$ is a frame for $\cH$, then there
%   exists a compact set $K\subseteq G$ (whose size depends only on the frame
%   bounds), such that $\cX \cap yK \neq  \emptyset $ for all $y\in G$. 
% \end{cor}

The comparison theorem is the starting point for genuine density
theorems. For instance, Landau's necessary conditions for set of
sampling are a simple corollary of Theorem~\ref{comp}. Likewise,
Ramanathan and Steger derived the necessary density of nonuniform
Gabor frames as an immediate corollary of a special case of 
Theorem~\ref{comp}. For more general groups and representations, one
may define some form of a Beurling density and deduce a density
theorem, but the results are not completely  satisfying.

\def\cprime{$'$} \def\cprime{$'$}

% \bibliographystyle{abbrv}
%  \bibliography{general,new}

\begin{thebibliography}{10}

\bibitem{BCHL06a}
R.~Balan, P.~G. Casazza, C.~Heil, and Z.~Landau.
\newblock Density, overcompleteness, and localization of frames. {I}. {T}heory.
\newblock {\em J. Fourier Anal. Appl.}, 12(2):105--143, 2006.

\bibitem{BCHL06}
R.~Balan, P.~G. Casazza, C.~Heil, and Z.~Landau.
\newblock Density, overcompleteness, and localization of frames. {II}. {G}abor
  systems.
\newblock {\em J. Fourier Anal. Appl.}, 12(3):309--344, 2006.

\bibitem{CDH99}
O.~Christensen, B.~Deng, and C.~Heil.
\newblock Density of {G}abor frames.
\newblock {\em Appl. Comput. Harmon. Anal.}, 7:292--304, 1999.

\bibitem{FNP07}
H.~Feichtinger, M.~Neuhauser, and M.~Piotrowski.
\newblock The homogeneous approximation property for coorbit spaces.
\newblock {\em In preparation}.

\bibitem{fg89jfa}
H.~G. Feichtinger and K.~Gr{\"o}chenig.
\newblock Banach spaces related to integrable group representations and their
  atomic decompositions. {I}.
\newblock {\em J. Functional Anal.}, 86(2):307--340, 1989.

\bibitem{FG07}
H.~F{\"u}hr and K.~Gr{\"o}chenig.
\newblock Sampling theorems on locally compact groups from oscillation
  estimates.
\newblock {\em Math. Z.}, 255(1):177--194, 2007.

\bibitem{GKT92}
K.~Gr{\"o}chenig, E.~Kaniuth, and K.~F. Taylor.
\newblock Compact open sets in duals and projections in {$L\sp 1$}-algebras of
  certain semi-direct product groups.
\newblock {\em Math. Proc. Cambridge Philos. Soc.}, 111(3):545--556, 1992.

\bibitem{GR96}
K.~Gr{\"o}chenig and H.~Razafinjatovo.
\newblock On {L}andau's necessary density conditions for sampling and
  interpolation of band-limited functions.
\newblock {\em J. London Math. Soc. (2)}, 54(3):557--565, 1996.

\bibitem{HK03}
C.~Heil and G.~Kutyniok.
\newblock Density of weighted wavelet frames.
\newblock {\em J. Geom. Anal.}, 13(3):479--493, 2003.

\bibitem{HK07}
C.~Heil and G.~Kutyniok.
\newblock The homogeneous approximation property for wavelet frames.
\newblock {\em J. Approx. Theory}, 147:28--46, 2007.

\bibitem{landau67}
H.~J. Landau.
\newblock Necessary density conditions for sampling and interpolation of
  certain entire functions.
\newblock {\em Acta Math.}, 117:37--52, 1967.

\bibitem{RS95}
J.~Ramanathan and T.~Steger.
\newblock Incompleteness of sparse coherent states.
\newblock {\em Appl. Comput. Harmon. Anal.}, 2(2):148--153, 1995.

\bibitem{sch98}
A.~P. Schuster.
\newblock The homogeneous approximation property in the {B}ergman space.
\newblock {\em Houston J. Math.}, 24(4):707--722, 1998.

\bibitem{Sun07}
W.~Sun.
\newblock Density of wavelet frames.
\newblock {\em Appl. Comput. Harmon. Anal.}, 22(2):264--272, 2007.

\bibitem{Sun07b}
W.~Sun.
\newblock Homogeneous approximation property for wavelet frames.
\newblock {\em Preprint}, 2007.

\bibitem{SunZ04}
W.~Sun and X.~Zhou.
\newblock Density of irregular wavelet frames.
\newblock {\em Proc. Amer. Math. Soc.}, 132(8):2377--2387 (electronic), 2004.

\end{thebibliography}

\end{document}